\documentclass[11pt,reqno]{amsart}
\usepackage[a4paper,margin=30mm]{geometry}
\usepackage{amsmath,amssymb,amsfonts,amsthm,mathtools}
\IfFileExists{mathrsfs.sty}{\usepackage{mathrsfs}}{\let\mathscr\mathcal}
\usepackage{microtype}
\usepackage[T1]{fontenc}
\usepackage{lmodern}
\usepackage{hyperref}
\hypersetup{colorlinks=true,linkcolor=blue,citecolor=blue,urlcolor=blue,pdftitle={Sharp Wasserstein Contractivity Exponents for the Porous Medium Equation in Several Dimensions},pdfauthor={Feizheng Huang and Zhaobo Liu}}
\usepackage[nameinlink,noabbrev,capitalise]{cleveref}
\allowdisplaybreaks
\newcommand{\RR}{\mathbb R}
\newcommand{\TT}{\mathbb T}
\newcommand{\PP}{\mathcal P}
\newcommand{\WW}{\mathsf W}
\newcommand{\Hs}{\mathscr Q}
\newcommand{\dd}{\mathrm d}
\newcommand{\tr}{\operatorname{tr}}
\newcommand{\supp}{\operatorname{supp}}
\newtheorem{theorem}{Theorem}[section]
\newtheorem{proposition}[theorem]{Proposition}
\newtheorem{lemma}[theorem]{Lemma}
\newtheorem{corollary}[theorem]{Corollary}
\theoremstyle{remark}
\newtheorem{remark}[theorem]{Remark}
\crefname{lemma}{Lemma}{Lemmas}
\crefname{proposition}{Proposition}{Propositions}
\crefname{corollary}{Corollary}{Corollaries}
\numberwithin{equation}{section}
\title[Sharp $\WW_p$ contractivity exponents for the PME]{Sharp Wasserstein Contractivity Exponents for the Porous Medium Equation in Several Dimensions}
\author{Feizheng Huang}
\author{Zhaobo Liu}
\usepackage{etoolbox}
\makeatletter
\patchcmd{\@setauthors}{\MakeUppercase{\authors}}{%
  {\normalfont\normalsize Feizheng Huang and Zhaobo Liu\par}
  \vspace{3pt}
  \begin{tabular}{@{}c@{\hspace{1.5em}}c@{}}
    \href{mailto:2024640043@mails.szu.edu.cn}{\color{black}\small\nolinkurl{2024640043@mails.szu.edu.cn}} &
    \href{mailto:liuzhaobo@szu.edu.cn}{\color{black}\small\nolinkurl{liuzhaobo@szu.edu.cn}}
  \end{tabular}\par
  \vspace{4pt}
  {\normalfont\small Institute for Advanced Study, Shenzhen University, Shenzhen, China\par}
}{}{\PackageError{authorblock}{Cannot place author contact details}{}}
\makeatother
\date{}
\begin{document}
\begin{abstract}
Let $S_t$ be the semigroup generated on $\RR^N$ by the porous medium equation
$\partial_tu=\Delta u^m$, where $N\ge2$ and $m>1$.
We determine exactly those $p\in[1,\infty]$ for which $S_t$ is nonexpansive
in the ordinary $p$-Wasserstein distance. For finite $p$ the classification
holds for all initial probability measures with finite $p$th moment.
The complete range is the closed interval
\[
\frac{2}{1+\vartheta_{m,N}}\le p\le\frac{2}{1-\vartheta_{m,N}},
\qquad
\vartheta_{m,N}\triangleq \left(\frac{1+N(m-1)}{m[1+(N-1)(m-1)]}\right)^{1/2}.
\]
For each $p$ outside this interval, including $1$ and $\infty$, there exist
two smooth radial probability densities, positive on the same ball and zero
outside, whose solutions have $\WW_p$ distance strictly greater than their
initial distance for all sufficiently small $t>0$.
\end{abstract}
\maketitle
\enlargethispage{2\baselineskip}

\section{Introduction and main result}
Let $m>1$ and $N\ge2$. We consider the Cauchy problem
\begin{equation}\label{eq:pme}
\partial_tu=\Delta(u^m),\qquad u(\cdot,0)=u_0\ge0,
\qquad x\in\RR^N.
\end{equation}
Write $S_t$ for the mass-preserving weak-solution semigroup and for its
action on nonnegative finite measures by initial traces
\cite[Thm.~4]{Pierre}. For $1\le p<\infty$, $\PP_p(\RR^N)$ denotes the
probability measures with finite $p$th moment. We use the ordinary distance
$\WW_p(\mu,\nu)=\inf_{\pi\in\Pi(\mu,\nu)}\||x-y|\|_{L^p(\pi)}$,
with the essential supremum when $p=\infty$. Here $\Pi(\mu,\nu)$ is the
set of couplings of $\mu$ and $\nu$. We identify a density $u$ with the
measure $u\,\dd x$.

In one space dimension all Wasserstein distances contract. In several
dimensions, $\WW_2$ nonexpansivity follows from the Wasserstein gradient-flow
structure \cite{Otto2001,OttoWest,BolleyCarrillo}. V\'azquez's focusing
example rules out even a uniform bound of the form $Ce^{\omega t}\WW_p(u_0,v_0)$
for $p=\infty$ and all sufficiently large finite $p$
\cite[Sec.~6, Theorem~4 and (6.1)]{Vazquez2005}.
Problem~1 of V\'azquez \cite[Sec.~7]{Vazquez2005} asks for the exact range
of exponents for which the distance between two solutions is nonincreasing.
For every $m>1$ and $N\ge2$, we show that this range is the closed interval
$[p_-,p_+]$, where
\begin{equation}\label{eq:theta}
\vartheta_{m,N}\triangleq 
\left(\frac{1+N(m-1)}{m[1+(N-1)(m-1)]}\right)^{1/2},\qquad
p_-\triangleq \frac2{1+\vartheta_{m,N}},\quad p_+\triangleq \frac2{1-\vartheta_{m,N}}.
\end{equation}
For $N\ge2$ and $m>1$ we have $0<\vartheta_{m,N}<1$, hence
$1<p_-<2<p_+<\infty$, and $1/p_-+1/p_+=1$.

\begin{theorem}[Sharp exponent classification]\label{thm:main}
Fix $m>1$ and $N\ge2$. For $1\le p<\infty$, the estimate
\begin{equation}\label{eq:global}
\WW_p(S_t\mu_0,S_t\nu_0)\le\WW_p(\mu_0,\nu_0)
\end{equation}
holds for every $\mu_0,\nu_0\in\PP_p(\RR^N)$ and every $t\ge0$ if
and only if $p\in[p_-,p_+]$, or equivalently
\begin{equation}\label{eq:criterion}
(N-1)(m-1)^2(p-2)^2\le4(p-1)[1+N(m-1)].
\end{equation}
For every $p\in[1,\infty]\setminus[p_-,p_+]$,
there are $u_0,v_0\in C_c^\infty(\RR^N)$, $u_0,v_0\ge0$ and
$\int u_0=\int v_0=1$, and a time $t_0>0$ such that
\begin{equation}\label{eq:strict}
\WW_p(S_tu_0,S_tv_0)>\WW_p(u_0,v_0)
\qquad\text{for every }0<t<t_0.
\end{equation}
For every such $p$ the two densities may be chosen radial, positive on the
same open ball and zero outside its closure.
\end{theorem}

By the semigroup law, \eqref{eq:global} for all initial pairs gives
nonincrease between any two times.

\begin{remark}\label{rem:range}
As $m\downarrow1$, the interval in \eqref{eq:theta} expands to
$(1,\infty)$, consistent with heat-flow contractivity. As $m\to\infty$ it
shrinks to $\{2\}$. Formally putting $N=1$ into \eqref{eq:criterion}
removes the restriction, consistent with the known one-dimensional result.
\end{remark}

The polynomial $P_p(b)$ in \eqref{eq:poly} links the two parts of the proof for finite $p$.
Inside the interval, its nonnegativity gives the sign of the matrix form
in the adjoint estimate. Outside the interval, a value of $b$ for which
$P_p(b)<0$ specifies a radial displacement proportional to $r^{-b}$ on
an annulus. For a suitable initial density, this displacement makes the
transport cost increase at the initial time.

\section{Contractivity in the sharp exponent range}\label{sec:global}
Write $a=m-1$, $k=N-1$ and, for $1<p<\infty$, $q=p/(p-1)$.
For $1\le p<\infty$ define the polynomial
\begin{equation}\label{eq:poly}
P_p(b)\triangleq m(p-1)b^2-akpb+k(1+ak).
\end{equation}
For real symmetric $N\times N$ matrices $H,K$ and a unit vector $e\in\RR^N$, write
\begin{equation}\label{eq:Q}
 \Hs(H,K)\triangleq\tr(KH^2)+a(\tr H)\tr(KH),\qquad
 K_q(e)\triangleq I+(q-2)e\otimes e .
\end{equation}

\begin{lemma}\label[lemma]{lem:algebra}
For $1<p<\infty$, \eqref{eq:criterion} is equivalent both to
$P_p\ge0$ on $\RR$ and to $\Hs(H,K_q(e))\ge0$ for every $H,e$.
When \eqref{eq:criterion} fails, $P_p$ takes a negative value at some $b>0$.
The last assertion also holds for $p=1$.
\end{lemma}
\begin{proof}
Take $e=e_1$, abbreviate $K_q(e_1)$ to $K_q$, and write $\alpha=H_{11}$, $\beta=k^{-1}\sum_{i>1}H_{ii}$.
Replacing $H$ by $\operatorname{diag}(\alpha,\beta,\ldots,\beta)$ preserves
$\tr H$ and $\tr(K_qH)$ and attains the lower bound
$\tr(K_qH^2)=\sum_{i,j}(K_q)_{ii}H_{ij}^2\ge(q-1)\alpha^2+k\beta^2$.
By homogeneity and continuity we may take $\beta=1$. Putting $\alpha=-(p-1)b$ gives
$\Hs(H,K_q)=k+(p-1)b^2+a(k-b)[k-(p-1)b]=P_p(b)$.
The discriminant of $P_p$ is
$k[ka^2(p-2)^2-4(p-1)(1+Na)]$, and its minimum occurs at
$b=akp/(2m(p-1))>0$. For $p=1$, $P_1(b)<0$ when $b>k+1/a$.
Finally $4(p-1)=p^2-(p-2)^2$ turns \eqref{eq:criterion}
into $|1-2/p|\le\vartheta_{m,N}$, giving \eqref{eq:theta}.
\end{proof}

We now use this matrix inequality to estimate solutions of the adjoint equation.
On the flat torus $\TT_L^N=(\RR/L\mathbb Z)^N$, let $\rho$ be a smooth
positive solution on $[0,T]$, and let $\phi$ solve the adjoint equation
with a smooth terminal value. Thus
\begin{equation}\label{eq:adjoint}
\rho_t=\Delta P,\qquad P=\rho^m,\qquad A=m\rho^{m-1},\qquad
\phi_t=-A\Delta\phi .
\end{equation}
The positive lower bound for $\rho$ ensures uniform parabolicity after
time reversal in the adjoint equation.

\begin{lemma}\label[lemma]{lem:cotangent}
Let $1<p<\infty$. If \eqref{eq:criterion} holds, then
$\int\rho(0)|\nabla\phi(0)|^q\le\int\rho(T)|\nabla\phi(T)|^q$.
\end{lemma}
\begin{proof}
For $h\in C^2(\RR^N)$ set $H=D^2\phi$ and $K=D^2h(\nabla\phi)$.
Using $\Delta h(\nabla\phi)=\tr(KH^2)+Dh(\nabla\phi)\cdot\nabla\Delta\phi$,
$P-\rho A=-aP$ and $\rho\nabla A=a\nabla P$ in \eqref{eq:adjoint} gives
\begin{align}
\frac{\dd}{\dd t}\int\rho h(\nabla\phi)
 &=\int P\tr(KH^2)-a\int Dh(\nabla\phi)\cdot\nabla(P\Delta\phi)\notag\\*
 &=\int P\,\Hs(H,K).\label{eq:Hidentity}
\end{align}
Take $h(z)=(|z|^2+\delta^2)^{q/2}/q$. Its Hessian is a positive scalar
times $(1-\theta)I+\theta K_q(e)$, where
$\theta=|z|^2/(|z|^2+\delta^2)$ and $e=z/|z|$ for $z\ne0$.
At zero $\theta=0$ and $e$ can be arbitrary. Since $\Hs(H,K)$ is linear
in $K$ and $I=K_2(e)$, \cref{lem:algebra} gives $\Hs(H,K)\ge0$
in \eqref{eq:Hidentity}. Integrate in time and let $\delta\downarrow0$ by
dominated convergence. This also covers equality in \eqref{eq:criterion}.
\end{proof}

We first prove contractivity on the torus and then pass to $\RR^N$.
\begin{proposition}\label[proposition]{prop:torus}
Under \eqref{eq:criterion}, the PME semigroup on $\TT_L^N$ is
$\WW_p$-nonexpansive on probability densities in $L^1$, for the geodesic
distance of the torus.
\end{proposition}
\begin{proof}
We interpolate between the initial measures along an optimal transport path
with parameter $s$ and smooth the entire path with the same heat kernel.
We then evolve each density for time $t$. The adjoint estimate controls
variation in $s$ along the evolved path and hence the distance between its endpoints.

Let $(X,Y)$ have an optimal coupling of the initial densities $u_0,v_0$.
Write $V$ for the representative of $Y-X$ in $[-L/2,L/2)^N$, so that
$D=(\mathbb E|V|^p)^{1/p}=\WW_p(u_0,v_0)$. Let $Z$ be independent of $(X,Y)$ with
density the torus heat kernel at time $\epsilon>0$. For $0\le s\le1$ set
\[
 X_s=X+sV+Z\pmod L,\qquad
 \rho_s(0,x)\,\dd x=\operatorname{Law}(X_s),\qquad \rho_s(t)=S_t\rho_s(0).
\]
Here $S_t$ denotes the PME semigroup on $\TT_L^N$.
Averaging translates of this positive smooth kernel gives initial
densities smooth in $(s,x)$ with a common positive lower bound, since $V$
is bounded. The maximum principle and uniform parabolicity give smooth
dependence on $s$ as in \cite[Sec.~4]{OttoWest}, with
$\partial_t\sigma_s=\Delta(m\rho_s^{m-1}\sigma_s)$ for
$\sigma_s=\partial_s\rho_s$. For each fixed smooth test function $f$,
differentiating $\mathbb E[f(X_s)]$ gives
$\int f\,\sigma_s(0)=\mathbb E[\nabla f(X_s)\cdot V]$.
Now fix $s$ and $t>0$, and solve \eqref{eq:adjoint} along $\rho_s$
with terminal value $\phi(t)=\psi\in C^\infty(\TT_L^N)$. The pairing $\int\phi(\tau)\sigma_s(\tau)$ is
constant for $0\le\tau\le t$. Taking $f=\phi(0)$ in the preceding
identity and applying H\"older's inequality and \cref{lem:cotangent} gives
\begin{equation}\label{eq:pair}
\left|\int\psi\,\sigma_s(t)\right|
 =\left|\mathbb E[\nabla\phi(0,X_s)\cdot V]\right|
 \le D\left(\int|\nabla\psi|^q\rho_s(t)\right)^{1/q}.
\end{equation}

At fixed $t$, take an optimal potential $\Phi$ in the torus dynamic dual
for the endpoints $\rho_0(t),\rho_1(t)$ and cost $|v|^p/p$
\cite[Prop.~2.2 and Eqs.~(2.2), (2.4)]{IshidaLavenant}.
It is Lipschitz in $(s,x)$ and satisfies
$\partial_s\Phi+|\nabla_x\Phi|^q/q\le0$ almost everywhere.
Write $\rho_s=\rho_s(t)$. Spatial mollification extends \eqref{eq:pair} to each Lipschitz slice
$\Phi(s,\cdot)$. With $g_s=\|\nabla_x\Phi(s,\cdot)\|_{L^q(\rho_s)}$,
the product rule and Young's inequality give, for almost every $s$,
\[
 \frac{\dd}{\dd s}\int\Phi(s,x)\rho_s(x)\,\dd x
 \le Dg_s-g_s^q/q\le D^p/p .
\]
By optimality, the endpoint difference is $\WW_p^p(\rho_0(t),\rho_1(t))/p$.
Integration over $0\le s\le1$ therefore gives $\WW_p(\rho_0(t),\rho_1(t))\le D$.
The initial endpoints are the heat convolutions of $u_0,v_0$ and converge
in $L^1$ as $\epsilon\downarrow0$. The PME is $L^1$-contractive, and on
the compact torus this convergence implies $\WW_p$ convergence.
\end{proof}

\begin{corollary}\label[corollary]{cor:measure}
For $p\in[p_-,p_+]$, \eqref{eq:global} holds on $\RR^N$ for all
$\mu_0,\nu_0\in\PP_p(\RR^N)$.
\end{corollary}
\begin{proof}
Write $B_r$ for the open ball of radius $r$ centered at the origin.
For smooth compactly supported data, finite propagation confines both flows
to some $B_R$ up to any fixed time $T$ \cite[Sec.~14.6.1]{VazquezBook}.
On a torus of side $L>4R$, periodization preserves the solutions by
uniqueness and all distances on $B_R\times B_R$. Thus \cref{prop:torus}
gives \eqref{eq:global}. Approximate arbitrary $\PP_p$ data in $\WW_p$ by
such densities. Pierre's theorem \cite[Thm.~4, p.~183]{Pierre} gives
$L^1$ convergence to the prescribed solutions at every $t>0$, and lower
semicontinuity of the extended transport cost gives \eqref{eq:global}.
Comparison with one compactly supported solution gives $S_t\mu_0\in\PP_p$.
\end{proof}

\section{Short-time growth outside the interval}\label{sec:nec}
By uniqueness, radial probability densities remain radial under the flow. We use the same letter for a radial function and its profile in $r=|x|$.
For radial profiles, primes denote differentiation with respect to $r$. The mass inside $B_r$ and the radius quantile are
\[
 M_u(r,t)\triangleq\int_{B_r}u(x,t)\,\dd x,\qquad
 R_u(s,t)\triangleq\inf\{r\ge0:M_u(r,t)\ge s\},\quad 0<s<1.
\]
For two such solutions and $1\le p\le\infty$,
\begin{equation}\label{eq:quantile}
 \WW_p(u(t),v(t))=\|R_v(\cdot,t)-R_u(\cdot,t)\|_{L^p(0,1)}.
\end{equation}
The common-direction quantile coupling is optimal for every finite $p$
\cite[Thm.~3.2 and Example~3.5]{CuestaEtAl}. Its maximal displacement
bounds $\WW_\infty$ from above. Taking $p\to\infty$ in
$\WW_p\le\WW_\infty$ gives the reverse bound.

\begingroup
\addtolength{\abovedisplayskip}{-2pt}
\addtolength{\belowdisplayskip}{-2pt}
\addtolength{\abovedisplayshortskip}{-2pt}
\addtolength{\belowdisplayshortskip}{-2pt}
Write $\omega_{N-1}=|\mathbb S^{N-1}|$. The radial flux gives
$M_{u,t}=\omega_{N-1}r^k(u^m)_r$ and $M_{u,r}=\omega_{N-1}r^ku$.
Enclose both initial quantile ranges of a compact mass interval in the
interior of a closed annulus with smooth positive data on a larger annulus.
Continuity at $t=0$ \cite[Sec.~7.5.1]{VazquezBook} keeps both densities
bounded away from zero on a neighborhood of the closed annulus for short
times. There $u_t=\nabla\cdot(m u^{m-1}\nabla u)$ is uniformly parabolic
with diffusivity smooth in $u$. Since the initial profile is smooth there,
local parabolic regularity up to the initial surface \cite{Lieberman}
gives $(u^m)_r(\cdot,t)\to(u_0^m)'$ uniformly on the closed annulus,
and likewise for $v$.
The $L^1$ trace gives uniform convergence of the ball masses, and the
positive lower bound for $M_{u,r}$ on the annulus gives uniform convergence
of the corresponding quantiles. Differentiating $M_u(R_u(s,t),t)=s$ yields
$\partial_t R_u=-M_{u,t}/M_{u,r}$ at $r=R_u(s,t)$. Thus the velocities
converge uniformly on the chosen mass interval as $t\downarrow0$, with
\begin{equation}\label{eq:quantile_velocity}
 \partial_tR_u(s,0+)=-\frac{(u_0^m)'(r)}{u_0(r)}
 =-\frac ma(u_0^a)'(r),\qquad
 r=R_u(s,0).
\end{equation}

\subsection{Finite exponents}
Fix $1\le p<\infty$ outside the interval and choose $b>0$ with $P_p(b)<0$.
Choose a smooth $d$ positive on $U_0=(3/4,5/4)$, zero off $U_0$, and equal
to $r^{-b}$ on a smaller interval $U_1$ about $1$. For small $\epsilon>0$,
$\tau(r)=r+\epsilon d(r)$ defines an increasing radial diffeomorphism
$T(x)=\tau(|x|)x/|x|$, with $T(0)=0$. Its Jacobian determinant is
$J(r)=\tau'(r)(\tau(r)/r)^k$ for $r>0$ and $J(0)=1$.
Let $u_0$ be a smooth radial probability density positive on $B_2$ and
zero outside, to be chosen below, and let $v_0=T_\#u_0$ be its pushforward under $T$. Then
$v_0(\tau(r))=u_0(r)/J(r)$ and $R_v(s,0)=\tau(R_u(s,0))$.

On the mass interval $I\triangleq(M_u(3/4,0),M_u(5/4,0))$, set
\[
 B(t)\triangleq\frac1p\int_I |R_v(s,t)-R_u(s,t)|^p\,\dd s .
\]
By \eqref{eq:quantile}, $B(t)\le \WW_p^p(u(t),v(t))/p$. Equality holds at $t=0$
because $\tau$ is the identity off $U_0$ and the initial quantiles agree outside $I$. It therefore suffices to prove $B'(0+)>0$.
Both initial
quantile ranges on $\overline I$ lie in $\overline{U_0}\subset(0,2)$.
The integrand's difference quotients are bounded by the velocity bounds and
local Lipschitz continuity of $|z|^p/p$. Positivity of the initial displacement
on $I$ gives pointwise right derivatives even at $p=1$, so dominated convergence applies.
With $w=d^{p-1}$ on $U_0$, \eqref{eq:quantile_velocity} and
$\dd s=\omega_{N-1}r^ku_0(r)\,\dd r$ give
\begin{equation}\label{eq:radgrowth}
 \frac{B'(0+)}{\omega_{N-1}\epsilon^{p-1}}=L_\epsilon(u_0^m),
\end{equation}
where, for each smooth profile $f$ on a neighborhood of $\overline{U_0}$,
\[
 L_\epsilon(f)\triangleq\int_{U_0}w\bigl[r^k f'-\tau^k(J^{-m}f)'\bigr]\,\dd r .
\]
Take a nonzero $g\in C_c^\infty(U_1)$, $g\ge0$, and extend it by zero. Integration by parts gives
$L_\epsilon(g)=\int_{U_1} g[J^{-m}(w\tau^k)'-(wr^k)']\,\dd r$ and $L_0(g)=0$.
Writing $\tau_\epsilon(r)=r+\epsilon d(r)$ and
$J_\epsilon(r)=\tau_\epsilon'(r)(\tau_\epsilon(r)/r)^k$, we obtain on $U_1$, where $d(r)=r^{-b}$,
\[
 \left.\partial_\epsilon J_\epsilon(r)\right|_{\epsilon=0}=(k-b)r^{-b-1}.
\]
In the derivative of $J^{-m}(w\tau^k)'$ at $\epsilon=0$, differentiating
$(w\tau^k)'$ and $J^{-m}$ gives, respectively, the two terms in the
coefficient of $r^{k-2-bp}$,
\[
 k(k-1-bp)-m(k-b)\bigl(k-b(p-1)\bigr)=-P_p(b).
\]
Hence
\[
 \left.\partial_\epsilon L_\epsilon(g)\right|_{\epsilon=0}
 =-P_p(b)\int_{U_1} g(r)\,r^{k-2-bp}\,\dd r>0 .
\]
Fix $\epsilon>0$ with $L_\epsilon(g)>0$. Take $u_0(x)$ proportional to
$(g(|x|)+\delta h(x)^m)^{1/m}$, with $h$ smooth, radial, positive in $B_2$ and
zero outside. The root is smooth near $\supp g$, where $h>0$, and
proportional to $h$ elsewhere. By linearity,
$L_\epsilon(g+\delta h^m)>0$ for small $\delta>0$.
Normalization preserves the sign, so $B'(0+)>0$ gives \eqref{eq:strict}.

\subsection{Maximal displacement}\label{sec:infinity}
Take $\tau(r)=r+R\chi(r)$, with $0\le\chi\le1$ smooth, supported in $U_0$
and equal to one near $1$. For small $R>0$, $\tau$ is increasing. Let $T$ be the radial map associated with this $\tau$.
Put $F(r)\triangleq 1-(r/(r+R))^{ka}$. Choose a smooth radial probability density $u_0$
positive in $B_2$ and zero outside, with $F(r)u_0^a(r)=Cr$ near $1$.
Since $F>0$ near $1$, this prescribes a smooth positive profile there,
and normalization only changes the constant $C>0$. Set $v_0=T_\#u_0$.
Near $1$, where $\tau(r)=r+R$, we have $v_0(r+R)=u_0(r)(r/(r+R))^k$.
For $s=M_u(r,0)$ there,
\[
 \partial_t(R_v-R_u)(s,0+)
 =\frac ma\,\frac{\dd}{\dd r}\bigl[u_0^a(r)-v_0^a(r+R)\bigr]
 =mC/a>0 .
\]
By \eqref{eq:quantile}, $\WW_\infty(u_0,v_0)=R$. Choose $[r_1,r_2]$ in this neighborhood with $r_1<r_2$.
Its image under $M_u(\cdot,0)$ has positive length since $u_0>0$.
Uniform velocity convergence on this mass interval gives
$R_v-R_u\ge R+mCt/(2a)$ for small $t>0$, proving \eqref{eq:strict}.

\begin{samepage}
\begin{remark}[Scaling]
Set $D_\lambda\mu=(x\mapsto x/\lambda)_\#\mu$ and $\gamma=2+Na$.
On a class closed under $D_\lambda$ for $\lambda\ge1$, apply a uniform
$Ce^{\omega t}$ bound to the dilated data at $t/\lambda^\gamma$.
Since $S_{t/\lambda^\gamma}D_\lambda=D_\lambda S_t$ and distances scale by
$\lambda^{-1}$, letting $\lambda\to\infty$ gives the factor $C$.
For each excluded $p$, \eqref{eq:strict} rules out $C=1$ for any finite
$\omega$ on every such class containing the corresponding counterexample pair.
\end{remark}
\end{samepage}

\par\endgroup
\enlargethispage{5\baselineskip}
\medskip
\indent\textbf{Acknowledgments.} The authors used OpenAI's ChatGPT and Codex to search the literature, explore revisions to the proofs, and prepare the manuscript. The authors assume responsibility for all content.

\indent A Lean 4 formalization verifying the main theorem and supporting results is available at\\
{\hypersetup{urlcolor=black}\urlstyle{same}\url{https://github.com/guohaoji21/sharp-wp-pme-lean}}.

\end{document}